\input amstex
\documentstyle{amsppt}
\magnification=\magstep 1
\catcode`\@=11
\def\nologo{\let\logo@=\relax}
\catcode`\@=\active
\nologo
\document
\vsize 7in
 
\topmatter

\title Socle degrees, Resolutions, and  Frobenius powers \endtitle 
 \leftheadtext{Kustin and Ulrich}
\rightheadtext{Socle degrees, Resolutions, and  Frobenius powers}
\author  Andrew R. Kustin\footnote{Supported in part by the National Security Agency.\hphantom{the National Science Foundation and XXX}} and Bernd  Ulrich\footnote{Supported in part by the National Science Foundation.\hphantom{and the National Security Agency XXX}}\endauthor
\address
Mathematics Department,
University of South Carolina,
Columbia, SC 29208 \endaddress
\email kustin\@math.sc.edu\endemail
\address Mathematics Department,
Purdue University,
West Lafayette, IN 47907
\endaddress
  \email ulrich\@math.purdue.edu  \endemail
\abstract 
We first  describe  a situation in which every graded Betti number in the tail of the resolution of $\frac RJ$ may be read from the socle degrees of $\frac RJ$. Then we apply the above result  to the ideals $J$ and $J^{[q]}$; and thereby 
describe  a situation in which the graded Betti numbers in the tail of the resolution of  $R/J^{[q]}$ are equal to the  graded Betti numbers in the tail of a shift 
 of the resolution of $R/J$.
\endabstract
\endtopmatter

Let\footnote""{\hskip-.17in 2000 {\it Mathematics Subject Classification.} Primary: 13D02, Secondary: 13A35.} $(R, \frak m)$ be a Noetherian graded algebra over a field of positive characteristic $p$, with irrelevant ideal $\frak m$.
Let $J$ be an $\frak m$-primary homogeneous ideal in $R$. Recall that if $q=p^e$, then the $e^{\text{th}}$ Frobenius power of $J$ is the ideal $J^{[q]}$ generated by all $j^q$ with $j\in J$. Recall, also,  that the socle of $\frac RJ$ is the ideal $\frac{(J \:\!\frak m )}{J}$ of $\frac RJ$.  The {\it socle degrees} of $\frac RJ$  are the degrees of any  
homogeneous basis for  the graded vector space $\operatorname{soc} \frac RJ$. 
The basic question is:

\proclaim{Question 0.1}
How do the socle degrees of $\frac{R}{J^{[q]}}$ vary with $q$? 
\endproclaim
The question of finding a linear bound for the top socle degree of $R/J^{[q]}$ has been considered by Brenner in \cite{1} from the point of view of finding inclusion-exclusion criteria for tight closure.
An answer to Question 0.1   would provide insight into the tight closure of $J$, and possibly a handle on Hilbert-Kunz functions.

We are particularly interested in how the socle degrees of Frobenius powers  encode homological information about the ideal $J$.
For example, the answer to   Question 0.1 is well-known in the case when $J$ has finite projective dimension: if the socle degrees of  $\frac RJ$ are $\{\sigma_i\mid 1\le i\le s\}$, then the socle degrees of $\frac R{J^{[q]}}$ are  $\{q\sigma_i-(q-1)a\mid 1\le i\le s\}$, where $a$ is the $a$-invariant of $R$. When $R$ is a complete intersection, the converse is established in \cite{6}. 

In the course of studying Question 0.1 we found some ideals $J$ for which it appeared that the tail of the minimal $R$-resolution of $R/J^{[p^{e}]}$ did not depend on $e$. In other words, there exist a complex $\Bbb  G_{\bullet}$, depending only on $J$, such that for each exponent $e$, there exists a twist $n_e$ and a complex 
$$\Bbb  F_{e,d-1}\to \Bbb  F_{e,d-2}\to \dots \to \Bbb  F_{e,0}$$ so that the minimal $R$-resolution of $R/J^{[p^e]}$ looks like 
$$\Bbb  G_{\bullet}(-n_e)\to \Bbb  F_{e,d-1}\to \Bbb  F_{e,d-2}\to \dots \to \Bbb  F_{e,0}, \tag 0.2$$where $d$ is the Krull dimension of $R$. 
Our examples were computer calculations made using Macaulay2. 
The Betti number charts showed that (0.2) might be possible. We performed row and column operations on the matrices in position $d+1$ to see that these matrices all presented the same $d^{\text{th}}$-syzygy, up to shift; hence, confirming (0.2) for the examples under consideration. 
To emphasize how special this phenomenon is, we point out  that we do not know the matrix in position $d$. Typically, the entries of this matrix have huge degrees (successive Frobenius powers grow very quickly); but magically,  the degrees of the entries in the  matrix in  position $d+1$ drop back down to the degrees that appeared in the resolution of $R/J$. Here is an example.
Let $P$ be the polynomial ring $\frac{\Bbb  Z}{(5)}[x,y,z]$, $f$ be the element $x^3+y^3+z^3$ of $P$,  $R$ be the hypersurface ring $P/(f)$, and  $J$ be the ideal $(x^5,y^5,z^5)$ of $R$. 
The  graded Betti numbers in the $R$-resolution of $R/J^{[p^e]}$ are:
$$\matrix 
e& \text{pos 0}& \text{pos 1}&\text{pos 2}&\text{pos 3}\\\noalign{\hrule}\\
0& 0\:\!1&5\:\!3&{8\:\!3\atop{9\:\!1}}&{9\:\!1\atop{10\:\!3}}\\\noalign{\hrule}\\
1& 0\:\!1&25\:\!3 &38\:\!3\atop{39\:\!1}&39\:\!1\atop{40\:\!3}\\\noalign{\hrule}\\
2& 0\:\!1&125\:\!3 &188\:\!3\atop{189\:\!1}&189\:\!1\atop{190\:\!3}\\\noalign{\hrule}\\
3& 0\:\!1&625\:\!3 &938\:\!3\atop{939\:\!1}&939\:\!1\atop{940\:\!3}\\\noalign{\hrule}\\
4& 0\:\!1&3125\:\!3 &4688\:\!3\atop{4689\:\!1}&4689\:\!1\atop{4690\:\!3.}\\\noalign{\hrule}\\
\endmatrix$$
The notation $5\:\!3$  means that the module in position $1$ in the minimal $R$-resolution of $R/J$ is $R(-5)^3$. Every entry in the matrix in position $2$ in the resolution of $R/J^{[p^4]}$ has degree $1563$ or $1564$; but every entry  of each matrix  in position $3$ is  linear or quadratic.   Row and column operations applied to the matrices calculated by Macaulay2 show that in each resolution the matrix in position $3$ may be taken to be 
$$\bmatrix 0&-x^2&-y^2&-2z\\x^2&0&-z^2&2y\\y^2&z^2&0&-2x\\2z&-2y&2x&0\endbmatrix.\tag 0.3$$

The present paper is a first attempt at identifying growth conditions on the socle degrees of $R/J^{[p^{e}]}$ which force  the tail of the minimal $R$-resolution of $R/J^{[p^{e}]}$ to be independent of $e$. 

The main result in this paper is 
Theorem 1.1 which describes a situation in which every graded Betti number in the tail of the resolution of $\frac RJ$ may be read from the socle degrees of $\frac RJ$. In Corollary 2.1, we apply  Theorem 1.1 to the ideals $J$ and $J^{[q]}$; and thereby 
describe  a situation in which the tail of the resolution of  $R/J^{[q]}$ is isomorphic to a shift of the tail of the resolution of $R/J$ as {\bf graded modules}. The relationship between the differentials in the two resolution tails remains a project for future study. 

Empirical evidence (see, for example, (0.3)) indicated that the map in position $3$ in the resolution of $R/J$ in Theorem 1.1 might be represented by an alternating matrix. In section 3, we prove that this does indeed happen. Section 4 consists of   a few examples and questions. In particular, we study the case where $R$ has small multiplicity and the socle of $R/J$ is pure. We find a situation where there is only one possibility for the dimension of the socle of $R/J$.  We also compare the $R$-resolutions of $R/J$ and the canonical module of ${R/J}$. 

The present paper has inspired further investigation of the phenomenon
\proclaim{Phenomenon 0.4} Sometimes, if the socle degrees of $R/J^{[q]}$ 
and $R/J$ are related ``correctly'', then the resolutions  of $R/J^{[q]}$ 
and $R/J$ share  the same tail. \endproclaim
One very interesting result along these lines is found in \cite{5}. Let
$R$ be the hypersurface ring $k[x,y]/(f)$, where $k$ is a field of characteristic $p$ and $f=x^n+y^n$, and let $J$ be the ideal $(x^N,y^N)R$. Conclusions (a), (b), and (c) are established in \cite{5}.
\proclaim\nofrills{(a)}\ The resolutions of $R/J$ and $R/J^{[q]}$ have a common tail, in the sense of {\rm(0.2)}, if and only if
$\operatorname{soc}(R/J^{[q]})$ and $[\operatorname{soc}(R/J)](-(q-1)N)$ are isomorphic as graded vector spaces. 
\endproclaim
\flushpar We emphasize that the common tail $\Bbb G_{\bullet}$ in the setting of \cite{5} is a complex {\bf with differential}. In other words, the first syzygy module $\operatorname{syz}_1^R(R/J^{[q]})$ is isomorphic to a shift of $\operatorname{syz}_1^R(R/J)$.
 Conclusion (a) confirms Phenomenon 0.4 -- including even the differential in the common tail  --  at least in the   set-up of \cite{5}.

\proclaim\nofrills{(b)}\ Once $N$, $n$, and $p$ are fixed, then,  there is a finite set of modules $\{M_i\}$ such that for each $q$ there exists an $i$ with   $\operatorname{syz}_1^R(R/J^{[q]})$  isomorphic to a shift of $M_i$. \endproclaim
 \flushpar Conclusion (b) is astounding. It leads to the natural question:
\proclaim{Question} What other rings $R$ and ideals $J$ have the property that  the set of syzygy modules $\{\operatorname{syz}_d^R(R/J^{[p^e]})\mid 0\le e\}$ is finite?\endproclaim
\flushpar In this Question, $d$ is the Krull dimension of the ring $R$ and two syzygy modules correspond to the same element of the finite set if one of the modules is isomorphic to a shift of the other module. 
\proclaim\nofrills{(c)}\ Furthermore, given $e$, there exist $n$, $p$, and $N$ such that 
$\operatorname{syz}_1^R(R/J^{[p^i]})$  are  distinct for  $0\le i\le e-1$, even after shifting; but
$\operatorname{syz}_1^R(R/J^{[p^k]})$   is isomorphic to a shift of $\operatorname{syz}_1^R(R/J^{[p^{k+e}]})$, for all $k$.
\endproclaim

\flushpar Conclusion (c) is a total surprise. We had been spoiled by the finite projective dimension case and the inability of the computer to compute more than a few successive Frobenius powers and  we had become used to having ``good behavior'' become visible fairly quickly. In fact, we now see that even when ``good behavior'' is guaranteed to occur, one might have to wait arbitrarily long until it is visible.

\bigskip We return to the present paper. Our situation is significantly more general than the situation of \cite{5}; and of course, our conclusions are not as explicit. In particular, we usually think of the common tail $\Bbb G$ as a graded module, not as a complex with differential.

We write $\omega_R$ for the graded canonical module of the graded ring $R$. See \cite{6} for information about canonical modules and $a$-invariants. If $m$ is a homogeneous element of a graded module $M$, then we write $|m|$ for the degree of  $M$. We use $\pmb s^n(\underline{\phantom{X}})$ to indicate that the degree of an element has been shifted by $n$. In other words, if $m$ is an element of the graded module $M$, and $n$ is an integer, then $\pmb s^{n}(m)$ is the element of $M(-n)$ which corresponds to $m$. In particular, $$|\pmb s^n(m)|=|m|+n.\tag 0.5$$So, $$m\in M_{|m|}\implies \pmb s^n(m)\in M(-n)_{|m|+n}.$$
 Let  $\mu(M)$ denote the minimal number of generators of the graded $R$-module $M$. We always use $(\underline{\phantom{X}})^*$ to mean the functor $\operatorname{Hom}_R(\underline{\phantom{X}},R)$; and we always use 
$\overline{\phantom {x}}$ to mean the functor $\underline{\phantom{X}} \otimes_PR$. 

If $I$ is a homogeneous Gorenstein ideal in a standard graded polynomial ring $P$, over a field,  then the last non-zero module in a minimal homogeneous resolution of $P/I$ by free $P$-modules has rank one and is equal to $P(-b)$ for some twist $b$. We refer to $b$ as the back twist  in the $P$-resolution of $P/I$ and we observe that the $a$-invariant of $P/I$ is equal to $b-\dim P$. 
 \heading Section 1. The main result. \endheading
In this section  $k$ is a field of arbitrary characteristic. The main result in this paper is 
Theorem 1.1 which describes a situation in which every graded Betti number in the tail of the resolution of $\frac RJ$ may be read from the socle degrees of $\frac RJ$. 
Remark 1.14 contains alternate versions of hypotheses (b) and (c).

\proclaim{Theorem 1.1} Let $P$ be a polynomial ring in three variables over a field $k$. Each variable has degree $1$. 
Let $f$ be a non-zero homogeneous element of $P$, $R$ be the hypersurface ring $R=P/(f)$, and $a$ be the $a$-invariant of $R$. Let $I$ be a homogeneous grade three Gorenstein ideal in $P$,  $b$ be the back twist in the $P$-resolution of $\frac PI$, and $J$ be the ideal $IR$. Let  $$\Bbb  F_{\bullet}:\quad \dots @>d_{4} >>\Bbb  F_{3}@>d_{3} >> \Bbb  F_{2}@>d_{2} >> \Bbb  F_{1}@> d_{1}>> R\to R/J\to 0\tag 1.2$$ be  the graded minimal $R$-resolution of $R/J$, and  $\{\sigma_{i}\mid 1\le i\le s\}$ be the socle degrees of $\frac RJ$.  Assume that 
\flushpar{\rm(a)} $I$ and $J$ have the same number of minimal generators,
\flushpar{\rm(b)} $\operatorname{rank} \Bbb  F_{2}=\dim_k\operatorname{soc}\frac RJ$, and
\flushpar{\rm(c)} $\sigma_{i}+\sigma_{j}\neq b+2a$ for any pair $(i,j)$.
\flushpar Then 
 \flushpar{\rm(A)} $\Bbb  F_{2}=\bigoplus\limits_{i=1}^{s}R(-(b+a-\sigma_{i}))$,  
\flushpar{\rm(B)} $ \Bbb  F_{3}=\bigoplus\limits_{i=1}^{s}R(-(\sigma_{i}+3))$, and
\flushpar{\rm(C)} $\Bbb  F_{i+2}=\Bbb  F_{i}(-|f|)$, for all $i\ge 2$.
\endproclaim
\demo{Proof} Let $Z=\operatorname{im} d_{2}$.

In the first part of the argument we identify a submodule $Z_1$ of $Z$ and prove that
$$\tsize \omega_{\frac RJ}(-b-a)\cong \frac{I\:\! f}{I}(-|f|)\cong\frac Z{Z_1}\tag 1.3$$ are isomorphic as graded $R$-modules.

There is no difficulty finishing the proof of (A)  once (1.3) has been established. Indeed, the definition of $Z$ yields that $\Bbb  F_{2}$ and $Z$ have the same generator degrees. We know (see, for example, 
\cite{6, Prop. 1.5}) that   $\mu(\omega_{\frac RJ})=\dim\operatorname{soc} \frac RJ$; furthermore, $$\text{the 
generator degrees of $\omega_{\frac RJ}$ are $\{-\sigma_{i}\mid 1\le i\le s\}$.}\tag 1.4$$ 
Use hypothesis (b) to see that
$$\mu(Z)= \mu(\Bbb  F_{2})= \dim\operatorname{soc} \tsize \frac RJ 
=\mu(\omega_{\frac RJ});$$
 and therefore, (1.3) shows that $\omega_{\frac RJ}(-b-a)$, $Z$, and $\Bbb  F_{2}$ all have the same generator degrees. Conclusion (A) now follows immediately.

We establish the left-most isomorphism of (1.3). Observe that the surjection $\frac PI\to \frac R{J}$ induces the equality
$$\tsize \omega_{\frac R{J}}=\operatorname{Hom}_{\frac PI}(\frac R{J},\omega_{\frac PI}).$$ 
The $a$-invariant of $\frac PI$ is $b-3$; so, the canonical module of $\frac PI$ is $\frac PI(b-3)$ and 
$$\tsize \omega_{\frac R{J}}=\operatorname{Hom}_{\frac PI}(\frac P{(I,f)},\frac PI(b-3))=\frac{I\:\! f}{I}(b-3). $$
Use the fact that $$a=|f|-3\tag 1.5$$ to see that the left-most isomorphism of (1.3) holds.

We define the submodule $Z_1$ of $Z$. Let $\operatorname{syz}_2^P(\frac PI)$ be the second syzygy of $\frac PI$ as a $P$-module and let 
$$\tsize \Bbb T_{\bullet}\:\quad 0\to \Bbb T_3@> t_3>> \Bbb T_2@> t_2>> \Bbb T_1@> t_1>> P\to \frac PI\to 0\tag 1.6$$ be the minimal $P$-resolution of $P/I$. Due to the  hypothesis $\mu(J)=\mu(I)$, we may choose $\Bbb T_1$ and $t_1$ so that $\Bbb T_1\otimes_PR=\Bbb  F_{1}$ and $t_1\otimes_PR=d_{1}$.
Apply $\underline{\phantom{X}}\otimes_PR$ to 
the exact sequence 
$$\tsize 0\to \operatorname{syz}_2^P(\frac PI)\to  \Bbb T_1@> t_1>> P\to \frac PI\to 0$$
and compare the resulting complex to the exact sequence that defines $Z$:
$$\eightpoint  \CD  @. \operatorname{syz}_2^P(\frac PI)\otimes_PR@>>>  \Bbb T_1\otimes_PR@> t_1\otimes_PR>> P\otimes_PR@>>> \frac PI\otimes_PR \\
@.@. @V = V\alpha_1V @V =VV @V =VV\\
0@>>> Z@>>>  \Bbb  F_{1}@> d_{1}>> R@>>> \frac R{J}.
\endCD\tag 1.7$$
The  bottom row is exact, so $\alpha_1$ induces a map  $\alpha_2\:\operatorname{syz}_2^P(\frac PI)\otimes_PR\to Z$. Let $Z_1=\alpha_2(\operatorname{syz}_2^P(\frac PI)\otimes_PR)$.

We establish the right-most isomorphism of (1.3). If $u\in (I\:\! f)(-|f|)$, then $uf=t_1\tau_1$ for some $\tau_1\in \Bbb T_1$; and therefore, $\tau_1\otimes R$ is in $\ker (t_1\otimes R)$ and $\alpha_1(\tau_1\otimes R)\in Z$. 
There is no difficulty in seeing that
$$h\:(I\:\! f)(-|f|)\to  \tsize \frac{Z}{Z_1},$$ 
given by
$$h(u)= \text{ the class of  $\alpha_1(\tau_1\otimes R)$,  where $uf=t_1\tau_1$ for some $\tau_1\in \Bbb T_1$},$$
is a well-defined homomorphism of graded $P$-modules. It is clear that $h$ is surjective and that $I(-|f|)$ is contained in the kernel of $h$. 

Now we show that the kernel of $h$ is contained in $I(-|f|)$. If   $u\in (I\:\! f)(-|f|)$ and $h(u)$ is zero in $\frac{Z}{Z_1}$, then there exists $\tau_1\in \Bbb T_1$ with $uf=t_1\tau_1$ and $\alpha_1(\tau_1\otimes R)=\alpha_1(\tau_1'\otimes R)$ in $Z\subseteq \Bbb  F_{1}$, for some $\tau_1'\in \operatorname{syz}_2^P(\frac PI)$. The hypothesis $\mu(I)=\mu(J)$ ensures that $\alpha_1$ is an isomorphism; hence, $(\tau_1-\tau_1')\otimes R$ is zero in $\Bbb T_1\otimes R$. In other words, there exists $\tau_1''$ in 
$\Bbb T_1$ with $\tau_1-\tau_1'=f\tau_1''$. Apply $t_1$ to see that $uf=ft_1(\tau_1'')$ in $P$. We know that $f$ is regular on $P$; so we conclude that $u=t_1(\tau_1'')\in I(-|f|)$. 
Both isomorphisms of (1.3) have been established, and the first part of the argument is complete.

The second syzygy module $Z$ is a maximal Cohen-Macaulay module over the two-dimensional ring $R$. A  straightforward calculation allows us to  to decompose   $$\Bbb  F_{3}@> d_{3}>> \Bbb  F_{2}@> d_{2}>> Z$$  into 
$$\Bbb  F_{3}@> {\bmatrix d_{3}'\\0\endbmatrix }>>\matrix \Bbb  F_{2}'\\\oplus\\  \Bbb  F_{2}''\endmatrix @> {\bmatrix d_{2}'&0\\0&d_{2}''\endbmatrix}>>  \matrix  Z'\\\oplus\\ Z''\endmatrix  $$
where $Z''$ is a free $R$-module, $Z'$ is a maximal Cohen-Macaulay $R$-module with no free summands, 
and  $d_{2}''$ is an isomorphism.
Eisenbud's ground breaking paper \cite{3} guarantees that the minimal 
resolution 
$$\dots @>d_{4} >>\Bbb  F_{3}@>d_{3}' >> \Bbb  F_{2}'@> d_{2}' >> Z'\to 0\tag 1.8$$ 
of $Z'$ by free $R$-modules is periodic of period $2$.  Furthermore, \cite{3} guarantees that the maps $d_{4}$ and $d_{3}'$ may be pulled back to $P$ to give a matrix factorization of $f$ times the identity matrix. 

In the second part of the argument, we show that 
$$\text{$Z^*(a)$ and $\omega_{R/J}$ have the same generator degrees}\tag 1.9$$
and
$$\text{$\Bbb  F_{3}^*(-3)$ and $Z^{\prime *}(a)$ have the same generator degrees.} \tag 1.10$$
Apply $\operatorname{Hom}_R(\underline{\phantom{X}},R(a))$ to the short exact sequence
$$0\to Z\to \Bbb  F_{1}\to J\to 0$$ to obtain the exact sequence
$$0\to J^*(a)\to \Bbb  F_{1}^*(a)\to Z^*(a)\to \operatorname{Ext}_R^1(J,R(a))\to 0.$$
Index shifting gives $\operatorname{Ext}_R^1(J,R(a))=\operatorname{Ext}^2_R(\frac R{J},R(a))$. The canonical module of $R$ is equal to  $\omega_R=R(a)$ and the canonical module of $\frac RJ$ is $\omega_{\frac R{J}}=\operatorname{Ext}_R^2(\frac R{J},\omega_R)$; see, for example, \cite{6, Prop. 1.2}.  We have produced a homogeneous degree zero surjection $$Z^*(a)\to \omega_{\frac R{J}}\to 0.\tag 1.11$$ Apply $\operatorname{Hom}_R(\underline{\phantom{X}},R)$ to (1.8) to see that
$Z^{\prime *}=\ker d_{3}^{\prime *}$. 
Extend the periodic resolution (1.8) to the left to obtain the homogeneous minimal resolution 
$$\dots@> d_{4}>> \Bbb  F_{3}@> d_{3}'>> \Bbb  F_{2}'@>d_{4}(|f|)>> \Bbb  F_{3}(|f|)@> d_{3}'(|f|)>> \Bbb  F_{2}'(|f|)@>>> Z'(|f|)\to 0\tag 1.12$$
The module $Z'$ is a maximal Cohen-Macaulay $R$-module; hence, $\operatorname{Ext}^i_R(Z',R)$ is zero for all positive $i$; and therefore, the dual of (1.12), which is
$$\eightpoint 0\to (Z'(|f|))^*\to (\Bbb  F_{2}'(|f|))^* @> (d_{3}'(|f|))^*>> (\Bbb  F_{3}(|f|))^*@>(d_{4}(|f|))^*>>(\Bbb F_2')^*@> d_{3}^{\prime *}>>\Bbb  F_{3}^*\to \dots \ ,$$
 is exact. We have produced an isomorphism of graded $R$-modules:
$$\frac {(\Bbb  F_{3}(|f|))^*}{\operatorname{im} (d_{3}'(|f|))^*} \cong Z^{\prime *}. $$
The resolution (1.8) is  minimal; hence, $\operatorname{im} (d_{3}'(|f|))^*$ is contained in $\frak m (\Bbb  F_{3}(|f|))^*$ and 
$(\Bbb  F_{3}(|f|))^*$ and $Z^{\prime *}$ have the same generator degrees. Use (1.5) to see that (1.10) holds.
 Furthermore, we also see that  
 $$\aligned &\tsize \mu(Z^*)= \mu(Z^{\prime *})+\mu(Z^{\prime\prime *})=\operatorname{rank} \Bbb  F_{3}^*+\operatorname{rank} \Bbb  F_{2}''=\operatorname{rank} \Bbb  F_{2}'+\operatorname{rank} \Bbb  F_{2}''\\ &\tsize=\operatorname{rank} \Bbb  F_{2}=\dim \operatorname{soc} \frac R{J}=\mu(\omega_{\frac RJ}).\endaligned\tag 1.13$$
Combine (1.13) and (1.11) to see that (1.9) holds.

 In the third part of the argument we prove that $Z''=0$. Suppose that $\frak z$ generates a free summand  of $Z$. Let $\zeta$ be a homogeneous element of $Z^*$ with $\zeta(\frak z)=1$. On the one hand, $\pmb s^{-a}(\zeta)$ is a minimal generator of $Z^*(a)$ of degree $-a-|\frak z|$. (The shift function is explained in (0.5).) So (1.9), together with  (1.4), yields that there exists $i$ with 
$$-a-|\frak z|=-\sigma_{i}.$$ On the other hand, in the first part of the argument, we already calculated that $$|\frak z|=b+a-\sigma_{j},$$ for some $j$. Hypothesis (c) guarantees that $\frak z$ does not exist. 
It follows that $Z''=0$, $Z=Z'$, and $\Bbb  F_{2}'=\Bbb  F_{2}$. Use (1.10), (1.9), and (1.4) to establish conclusion (B) and (1.8) to establish conclusion (C). 
\qed
\enddemo

\remark{Remark 1.14}The isomorphism 
$$\tsize \omega_{\frac RJ}(-b-a)\cong\frac Z{\operatorname{im}(\operatorname{syz}_2^P(\frac PI)\otimes_PR)}$$of 
(1.3) would continue to hold even if the hypothesis
$\operatorname{rank} \Bbb  F_{2}=\dim_k\operatorname{soc}\frac RJ$
had not been imposed. In other words, in the context of Theorem 1.1, one always has
$$\tsize \dim \operatorname{soc} \frac RJ\le \operatorname{rank} \Bbb  F_{2},$$ and equality holds if and only if every  $P$-syzygy of $I$ is inside the maximal ideal  times the module of  $R$-syzygies of $J$.  

We use 
hypothesis (c) to prove

\proclaim\nofrills{(c$'$)}\ the second syzygy of the $R$-module $R/J$ does not have any free summands.\endproclaim

\flushpar Theorem 1.1 remains valid if one replaces hypothesis (c) with hypothesis (c$'$). However, there are two advantages to (c) over (c$'$). First, (c) is given in terms of the data $\{\sigma_i\}$, $a$, and $b$ of Theorem 1.1. Second,   in the proof of Corollary 2.1  we apply Theorem 1.1 to both $R/J$ and $R/J^{[q]}$. If we assume that (c) holds for $R/J$, then the corresponding formula automatically holds for $R/J^{[q]}$, without making any further assumption. 
\endremark

 \heading  Section 2. The application to Frobenius powers. \endheading

The present paper was motivated by the observation of Phenomenon 0.4:   sometimes, if the socle degrees of $R/J^{[q]}$ 
and $R/J$ are related ``correctly'', then the resolutions  of $R/J^{[q]}$ 
and $R/J$ share the same tail,   
after a shift. In Corollary 2.1, we apply   Theorem 1.1  twice and obtain a situation where the tail of the resolution of $R/J^{[q]}$ is a shift of the tail of the resolution of $R/J$,   as {\bf  graded modules}.

\proclaim{Corollary 2.1} Let $P$ be a polynomial ring in three variables over a field $k$ of positive characteristic $p$. Each variable has degree $1$. 
Let $f$ be a non-zero  homogeneous element of $P$, $R$ be the hypersurface ring $R=P/(f)$,  and $a$ be the $a$-invariant of $R$.  Let $I$ be a homogeneous grade three Gorenstein ideal in $P$,  $b$ be the back twist in the $P$-resolution of $\frac PI$, $J$ be the ideal $IR$, $\{\sigma_{i}\mid 1\le i\le s\}$ be the socle degrees of $\frac RJ$, and $\Bbb  F_{\bullet}$, as given in {\rm(1.2)}, be the graded minimal $R$-resolution of $R/J$. Let $e$ be a fixed exponent, $q=p^e$,  and $$\Bbb  F_{e,\bullet}:\quad \dots @>d_{e,4} >>\Bbb  F_{e,3}@>d_{e,3} >> \Bbb  F_{e,2}@>d_{e,2} >> \Bbb  F_{e,1}@> d_{e,1}>> R\to R/J^{[q]}\to 0$$ be  the graded minimal $R$-resolution of $R/J^{[q]}$. 
Assume that \flushpar{\rm (a)} $I$, $J$, and $J^{[q]}$ have the same number of minimal generators, 
    \flushpar{\rm (b)} $\operatorname{rank} \Bbb  F_{2}=\dim_k\operatorname{soc}\frac RJ$, and
\flushpar{\rm (c)} $\sigma_{i}+\sigma_{j}\neq b+2a$ for any pair $(i,j)$.
\flushpar 
If
$$\tsize \operatorname{soc} \frac R{J^{[q]}} \quad\text{and}\quad  (\operatorname{soc} \frac RJ)\left(-\frac {b(q-1)}2\right)\tag 2.2$$ are isomorphic as graded vector spaces, then 
 $$\tsize \Bbb  F_{e,i}\quad\text{and}\quad  \Bbb  F_{i}\left(-\frac {b(q-1)}2\right)$$ are isomorphic as  graded modules
for all integers $i\ge 2$.
\endproclaim

\demo{Proof}We first show  that  $\operatorname{rank} \Bbb  F_{e,2}=\dim_k\operatorname{soc}\frac R{J^{[q]}}$. In light of Remark 1.14, this amounts to showing that every $P$-syzygy of $I^{[q]}$ is inside the maximal ideal  times the module of  $R$-syzygies of $J^{[q]}$. This does happen because the hypothesis tells us that
 every   $P$-syzygy of $I$ is inside the maximal ideal  times the module of  $R$-syzygies of $J$, and the  Frobenius homomorphism is flat on $P$-modules so every $P$-syzygy of $I^{[q]}$ is the $q^{\text{th}}$-Frobenius power of a  $P$-syzygy of $I$. 

Notice also that the analogue of hypothesis (c) holds for the ideal $J^{[q]}$ of $R$. One consequence of hypothesis (2.2) is that the socles of $R/J$ and $R/J^{[q]}$ have the same dimension.  Let $\{\sigma_{e,i}\mid 1\le i\le s\}$ be the socle degrees of $R/J^{[q]}$. Hypothesis (2.2) yields that 
$\sigma_{e,i}+\sigma_{e,j}=\sigma_{i}+\sigma_{j}+b(q-1)$. We know that the
back twist in the $P$-resolution of $\frac P{I^{[q]}}$ is
$$b_e=qb;\tag 2.3$$ therefore, hypothesis (c) guarantees that
$\sigma_{e,i}+\sigma_{e,j}\neq b_e+2a$ for any pair $(i,j)$.

We apply Theorem 1.1 to $I$ and to $I^{[q]}$.  It follows from (2.3) that
$$\tsize b_e-\frac {b(q-1)}2=b+\frac {b(q-1)}2;$$ and therefore
$$\tsize \Bbb  F_{e,2}=\bigoplus\limits_{i=1}^{s}R(-(b_e+a-\sigma_{i}-\frac {b(q-1)}2))=\Bbb  F_{2}\left(-\frac {b(q-1)}2\right).$$ The calculation $\Bbb  F_{e,3}=\Bbb  F_{3}\left(-\frac {b(q-1)}2\right)$ is even easier. Both resolutions $\Bbb  F_{e,\bullet}$
and $\Bbb  F_{\bullet}$ are eventually periodic with 
$\Bbb  F_{e,i+2}=\Bbb  F_{e,i}(-|f|)$  and $\Bbb  F_{i+2}=\Bbb  F_{i}(-|f|)$ for $i\ge 2$. \qed

\enddemo
\remark{Remark} If the polynomial $f$ of Corollary 2.1 is irreducible and   all of the minimal generators of $I$ have the same degree, then the hypothesis $\mu(J^{[q]})=\mu(J)$ automatically holds. Indeed, inflation of the base field $k\to K$ gives rise to faithfully flat extensions $P\to P\otimes_kK$ and $R\to R\otimes_kK$. Consequently,  we may assume  that $k$ is a perfect field. If  $g$ is a minimal generator of $I^{[q]}$, then the hypothesis that all of the minimal generators of $I$ have the same degree ensures that     $g=\sum\limits_{i=1}^n\alpha_ig_i^{[q]}$, where $(g_1,\dots,g_n)$ is a minimal generating set for $I$, and each $\alpha_i$ is in $k$.  The field $k$ is perfect; so each $\alpha_i$ has a $q^{\text{th}}$-root and $g=g_0^q$, for some minimal generator  $g_0$ of $I$. If $g$ were in $(f)$, then $g_0$ would also be in $(f)$, since $(f)$ is a prime ideal. \endremark

\remark{Remark} The hypothesis of Corollary 2.1 is far from arbitrary; that is,  the only possible number $n$ with
$$\Bbb  F_{e,3}=\Bbb  F_{3}(-n)\quad 	\text{and}\quad \Bbb  F_{e,2}=\Bbb  F_{2}(-n)\tag 2.4$$ is $n=\frac{b(q-1)}2$. Indeed, if (2.4) occurs, then Theorem 1.1 shows that 
$$-\sigma_{e,i}-3=-\sigma_{i}-3-n \quad 	\text{and}\quad -b_e-a+\sigma_{e,i}=-b-a+\sigma_{i}-n,\tag 2.5$$
where $\{\sigma_{e,i}\}$ and $\{\sigma_{i}\}$ are the socle degrees of $R/J^{[q]}$ and $R/J$, respectively.  One may solve (2.5) to see that
$n$ must equal  $\frac{b(q-1)}2$ and $\operatorname{soc} \frac R{J^{[q]}}$ must equal $\operatorname{soc} \frac RJ(-n)$.
\endremark

\heading  Section 3. Alternating maps in the resolution of $R/J$. \endheading

This section is a continuation of section one; the field $k$ may have any characteristic.   Corollary 3.1, which is the main result in the section,
establishes that, for any ideal as in Theorem 1.1, the maps 
in the tail of the resolution are alternating. The proof appears at the end of the section. 

\proclaim{Corollary 3.1}If the notation and hypotheses of Theorem 1.1 are in effect, then there exist homogeneous alternating $s\times s$ matrices $\Phi$ and $\Psi$ with entries in $P$ such that
 $\Phi\Psi=fI=\Psi\Phi$ and 
$$\dots @>\psi>>\Bbb  F_{3}(-|f|)@> \varphi>> \Bbb  F_{2}(-|f|) @> \psi>> \Bbb  F_{3} @> \varphi>> \Bbb  F_{2}@> d_{2}>> \Bbb  F_{1}@> d_{1}>> R\tag 3.2$$ is the minimal homogeneous resolution of $R/J$ by free $R$-modules, where $\varphi=\overline{\Phi}$ and $\psi=\overline{\Psi}$.\endproclaim 

The determinant of $\Phi$ cannot be zero; so one consequence of Corollary 3.1 is that $s$ must be even. 

There are three steps in the proof of Corollary 3.1. 
The proof of Theorem 1.1 depended on very careful analysis of $Z$ and  $Z^*$, where $Z$ is the second syzygy module of $R/J$. Now that the proof of Theorem 1.1 is complete, we are able to prove that $Z^*$ is isomorphic to a shift of $Z$, and this is the first step in the proof of Corollary 3.1. The isomorphism that we produce in Corollary 3.3 appears to be fairly abstract; however, at least part of it is induced by the alternating map in the Buchsbaum-Eisenbud resolution of the grade three Gorenstein ideal $I$. Furthermore,   there exist a nonzerodivisor  $c$ in $R$  with $cZ$ contained in the part of $Z$ on which the isomorphism  is known to be alternating. This calculation appears in Proposition 3.10. The final step uses ideas from \cite{4}. 

\proclaim{Corollary 3.3} Retain the notation and hypotheses of Theorem 1.1 and its proof. Then the graded $R$-modules $$ \tsize Z(b),\quad (\operatorname{syz}_2^P(\frac PI)\otimes_PR)^*, \quad \text{and}\quad Z^*$$ are isomorphic.\endproclaim

\demo{Proof} Recall the $P$-resolution $\Bbb T_{\bullet}$ of the codimension three  Gorenstein ring $P/I$ which is given in (1.6). The module $\Bbb T_3$ is equal to $P(-b)$ and the entries of the matrix  $t_3$ generate the ideal $I$. 
The $P$-resolution of $\operatorname{syz}_2^P(P/I)$ is 
$$0\to  \Bbb T_3\to  \Bbb T_2\to \operatorname{syz}_2^P(P/I)\to 0,$$ and 
$$\operatorname{Tor}_1^P(\operatorname{syz}_2^P(P/I),R)=\operatorname{Tor}_3^P(P/I,R)=0;$$ hence, the $R$-resolution of $\operatorname{syz}_2^P(P/I)\otimes_PR$ is 
$$0\to  \Bbb T_3\otimes_PR\to  \Bbb T_2\otimes_PR\to \operatorname{syz}_2^P(P/I)\otimes_PR\to 0.\tag 3.4$$ 
Apply $(\underline{\phantom{X}})^*=\operatorname{Hom}_R(\underline{\phantom{X}},R)$ to the short exact sequence 
(3.4) to see that 
$$\tsize 0\to (\operatorname{syz}_2^P(\frac PI)\otimes_PR)^*\to (\Bbb T_2\otimes_PR)^*@>(t_3\otimes 1)^*>> (\Bbb T_3\otimes_PR)^*\to \tsize \frac RJ(b)\to 0$$ is exact. We have established that
$$\tsize Z(b)\cong (\operatorname{syz}_2^P(\frac PI)\otimes_PR)^*.\tag 3.5$$ 

The isomorphism of (1.3) may be reconfigured as a short exact sequence 
$$0\to \operatorname{syz}_2^P(P/I)\otimes_PR @> \alpha_2>> Z\to \omega_{R/J}(-b-a)\to 0.\tag 3.6$$
(The map $\operatorname{syz}_2^P(P/I)\otimes_PR\to  \Bbb T_1\otimes_PR$ of (1.7) is injective because $\operatorname{Tor}^P_1(I,R)= \operatorname{Tor}^P_2(P/I,R)=0$.) The grade of the annihilator of the $R$-module $\omega_{R/J}$ is two; hence, $\operatorname{Ext}^i_R(\omega_{R/J},R)=0$ for $0\le i\le 1$. 
Apply  $(\underline{\phantom{X}})^*$ to the exact sequence (3.6) to conclude that $$\tsize \alpha_2^*\:Z^*\to (\operatorname{syz}_2^P(\frac PI)\otimes_PR)^* \tag 3.7$$ is an isomorphism.  \qed \enddemo

By carefully analyzing the isomorphisms (3.5) and (3.7), we are now able to learn much more about the isomorphism $Z^*\to Z(b)$ from Corollary  3.3. Keep the minimal $P$-resolution $({\Bbb T},t)$
of $P/I$ from (1.6). Fix an orientation isomorphism $[\underline{\phantom{X}}]\:  {\Bbb T}_3\to P(-b)$. Buchsbaum and Eisenbud \cite{2} proved that $ \Bbb T$ has the structure of a DG$\Gamma$-algebra.  In particular, the maps $\eta_i\:{\Bbb T}_i\to \operatorname{Hom}_P({\Bbb T}_{3-i},P(-b))$, which are given by
$$\eta_i(\tau_i)= [\tau_i\cdot \underline{\phantom{X}}],\tag 3.8$$for all $\tau_i\in \Bbb T_i$, give rise to an isomorphism of complexes $\Bbb T_{\bullet}\cong \operatorname{Hom}_{P}(\Bbb T_{\bullet},P(-b))$. 
Furthermore, if $\tau_2\in {\Bbb T}_2$, then $$t_2(\tau_2)\cdot \tau_2=t_4(\tau_2^{(2)})=0.\tag 3.9$$ Recall that $\overline{\phantom {X}}$ is the functor $\underline{\phantom{X}} \otimes _PR$. The DG$\Gamma$-structure on $\Bbb T$ induces a DG$\Gamma$-structure on $\overline{\Bbb T}$.  

In this discussion we take $Z$ to be $\frac {\Bbb F_{2}}{\operatorname{im} d_{3}}$.  
Every element of $Z$ is equal to $\operatorname{\frak q}(y)$, for some element  $y$  of $\Bbb F_{2}$, where $\operatorname{\frak q}\:\Bbb F_{2} \to \frac {\Bbb F_{2}}{\operatorname{im} d_{3}}=Z$ is the natural quotient map.  The shift function $\pmb s^n$ is explained in (0.5).

\proclaim{Proposition 3.10} Let $\ell\:Z^*\to Z(b)$ be one of the isomorphisms  from Corollary  3.3.  \item{\rm (1)}One may choose $\ell$ so that for all  $\zeta\in Z^*$ and $y\in\Bbb F_{2}$, the following statements are equivalent
\itemitem{\rm (a)} $\pmb s^{b}(\ell( \zeta))=\operatorname{\frak q}(y)$ in $Z$
\itemitem{\rm (b)}  $\pmb s^{-b}([d_{2}(y)\cdot \bar\tau_2])=\zeta(\operatorname{\frak q}(y'))$ in $R$ for all pairs $(\bar\tau_2,y')\in \overline{ \Bbb T}_2\times \Bbb F_{2}$ with $\bar t_2(\bar\tau_2)=d_{2}y'$ in $\overline {\Bbb T}_1=\Bbb F_{1}$.
\item{\rm (2)}The map $\ell$ from {\rm(1)} is an alternating map in the sense that $\zeta(\pmb s^{b}\ell(\zeta))$ is equal to zero for all $\zeta$ in $Z^*$.\endproclaim
\demo{Proof} We first prove (1). The map $\ell$ is the composition of isomorphisms
$$Z^*@>\alpha_2^*>> (\operatorname{syz}_2^P(\tsize\frac PI)\otimes_PR)^*@>\gamma >> Z(b),$$where $\alpha_2$ is defined in (1.7) and $\gamma$ is one the isomorphisms of (3.5). 
The map $\alpha_2^*$ is shown to be an isomorphism in (3.7). Eventually, we will make   a particular choice for $\gamma$.  
Fix $\zeta\in Z^*$. Every element of 
$Z(b)$ has the form $\pmb s^{-b}(\operatorname{\frak q}(y))$ for some $y\in \Bbb F_2$. Fix one such $y$. We compare the elements $\alpha_2^*(\zeta)$ and $\gamma^{-1}(\pmb s^{-b}(\operatorname{\frak q}(y)))$ of $(\operatorname{syz}_2^P(\tsize\frac PI)\otimes_PR)^*$.
For $R$-modules $M$ and $N$, we let 
$${<}\phantom{X},\phantom{X}{>}\:\operatorname{Hom}_R(M,N)\otimes M \to N$$ represent the evaluation map.  We write
 $\operatorname{syz}_2^P(\frac PI)$ as $\frac{\Bbb T_2}{\operatorname{im} t_3}$.  Let    
$$\tsize \operatorname{\frak q}'\:\Bbb T_2\to \frac{\Bbb T_2}{\operatorname{im} t_3}= \operatorname{syz}_2^P(\frac PI)$$ be the natural quotient map. A typical element of $\operatorname{syz}_2^P(\tsize \frac PI)\otimes R$ has the form $\overline {\frak q'(\tau_2)}$ for some $\tau_2\in \Bbb T_2$. We see that
$$\split \text{condition (a) holds} &{} \iff \ell(\zeta)=\pmb s^{-b}(\operatorname{\frak q}(y))\\&{}\iff 
\alpha_2^*(\zeta)= \gamma^{-1}(\pmb s^{-b}(\operatorname{\frak q}(y)))\\&\iff
{<}\alpha_2^*(\zeta),\overline {\frak q'(\tau_2)}{>}={<}\gamma^{-1}(\pmb s^{-b}\operatorname{\frak q}(y)),\overline {\frak q'(\tau_2)}{>},\\&\phantom{{} \iff{}} \text{for all $\tau_2\in \Bbb T_2$.}
\endsplit$$We next compute ${<}\alpha_2^*(\zeta),\overline {\frak q'(\tau_2)}{>}$ and ${<}\gamma^{-1}(\pmb s^{-b}\operatorname{\frak q}(y)),\overline {\frak q'(\tau_2)}{>}$.
 
If $\tau_2\in \Bbb T_2$, then $t_2(\tau_2)$ is in the kernel of $t_1$ and $\overline{t_2(\tau_2)}$ is in the kernel of  $ \bar t_1=d_{1}$, see (1.7); hence, there is an element $y'$ in $\Bbb  F_{2}$ with $d_{2}(y')=\overline{t_2(\tau_2)}$.  It follows that $\alpha_2\left(\overline{\operatorname{\frak q}'(\tau_2)}\right)=\operatorname{\frak q}(y')$; and therefore, 
$${<}\alpha_2^*(\zeta),\overline {\frak q'(\tau_2)}{>}=\zeta(\operatorname{\frak q}(y')).$$

 The argument of (3.5) shows that there exists an isomorphism  $\xi$ for which the diagram 
$$\eightpoint\CD Z@>d_2>>\Bbb F_1=\overline{\Bbb T}_1 @> \bar t_1 >> \overline{\Bbb T}_0\\@V\xi VV @V \eta_1 VV @VVV\\
\operatorname{Hom}_R\left(\operatorname{syz}_2^P(\frac PI)\otimes R,R(-b)\right)@>{\overline{\operatorname{\frak q}'}}^*>> \operatorname{Hom}_R(\overline {\Bbb T}_2,R(-b))@>\bar t_3^*>>\operatorname{Hom}_R(\overline {\Bbb T}_3,R(-b))
\endCD\tag{3.11}$$commutes, where $\eta_1$ is defined in (3.8). 
Twist the isomorphism $$\tsize\xi\: Z\to \operatorname{Hom}_R\left(\operatorname{syz}_2^P(\frac PI)\otimes R,R(-b)\right)= (\operatorname{syz}_2^P(\frac PI)\otimes R)^*(-b)$$ by $b$ to obtain the isomorphism 
$$\tsize\xi(b)\: Z(b)\to(\operatorname{syz}_2^P(\frac PI)\otimes R)^*.$$ At this point we choose $\gamma$ to be the inverse of $\xi(b)$. Therefore,
$$\gamma^{-1}(\pmb s^{-b}\operatorname{\frak q}(y))=(\xi(b))(\pmb s^{-b}\operatorname{\frak q}(y))=\pmb s^{-b}\left( \xi (\operatorname{\frak q}(y))\right).$$
It follows that
$${<}\gamma^{-1}(\pmb s^{-b}\operatorname{\frak q}(y)),\overline {\frak q'(\tau_2)}{>}
={<}\pmb s^{-b}\left( \xi (\operatorname{\frak q}(y))\right),\overline {\frak q'(\tau_2)}{>}
=\pmb s^{-b}({<}\xi (\operatorname{\frak q}(y)),\overline {\frak q'(\tau_2)}{>}).
$$ Follow the commutative square in (3.11) which defines $\xi$ to see that 
$${<} \xi (\operatorname{\frak q}(y)),\overline {\frak q'(\tau_2)}{>}=[d_{2}(y)\cdot \bar \tau_2];$$thus,
$${<}\gamma^{-1}(\pmb s^{-b}\operatorname{\frak q}(y)),\overline {\frak q'(\tau_2)}{>}=
\pmb s^{-b}([d_{2}(y)\cdot \bar \tau_2]),$$ and the proof of (1) is complete.

The proof of (2) has two parts. We first show that $\zeta(\pmb s^{b}\ell(\zeta))=0$ whenever $\zeta$ is in $Z^*$ and  $\pmb s^{b}\ell(\zeta)$ is an element of the submodule $Z_1$ (see (1.7)) of $Z$; then we show that there is a nonzerodivisor $c$ in $R$ with $cZ\subseteq Z_1$. Let $\zeta\in Z^*$ and let $y$ be an element of $\Bbb  F_{2}$ with $\pmb s^{b}\ell(\zeta)=\operatorname{\frak q}(y)$. Suppose first that there exists an element $\bar \tau_2$ of $\overline{\Bbb T}_2$ with $\bar t_2(\bar \tau_2)=d_{2}(y)$. In this case, we use (1) to see that
$$\zeta(\pmb s^{b}\ell(\zeta))=\zeta(\operatorname{\frak q} (y))=\pmb s^{-b}([d_{2}(y)\cdot \bar\tau_2])=\pmb s^{-b}([\bar t_2(\bar \tau_2)\cdot \bar \tau_2])=0.$$The final equality follows from  (3.9).  
We turn to the second part of the argument. The grade two ideal $J=IR$ is not contained in any associated prime ideal of the hypersurface ring $R$. We show that if $c$ is any element of $I$ and $y$ is any element of $\Bbb  F_{2}$, then $cd_{2}(y)$ is in the image of $\bar t_2$. Let $\tau_1$ be any lifting of the element $d_{2}(y)$ from $\Bbb  F_{2}=\Bbb T_1\otimes_PR$ back to  $\Bbb T_1$. The fact that  $d_{2}(y)$ is in the kernel of $d_{1}=\bar t_1$ ensures that there is an element $u$ in $P$ with $t_1(\tau_1)=uf$. There is an element $\tau_1'$ in $\Bbb T_1$ with $t_1(\tau_1')=c$. We see that
$$ct_1(\tau_1)=t_1(\tau_1')uf;$$ and therefore,
$c\tau_1-uf\tau_1'$ is in $\ker t_1=\operatorname{im} t_2$. Apply $\overline{\phantom {X}}=\underline{\phantom{X}} \otimes _PR$ to see that
$cd_{2}(y)= c\bar \tau_1\in \operatorname{im} \bar t_2$. 
 \qed 
\enddemo

\demo{Proof of Corollary 3.1}
We saw at the end of the proof of Theorem 1.1 that $Z$ is a maximal Cohen-Macaulay $R$-module with no free summands and therefore  $Z$ has a periodic resolution of period two which is induced by a matrix factorization of $f$. Indeed, the proof of \cite{3, Thm\.~6.1} shows that there exist matrices $D_3$ and $D_4$ over $P$ so that\roster
\item"{(1)}"
 $$  D_3\:\bigoplus\limits_{i=1}^sP(-\sigma_i-3)\to \bigoplus_{i=1}^sP(\sigma_i-b-a)
\quad \text{and}\tag{3.12}$$ $$ D_4\:\bigoplus\limits_{i=1}^sP(\sigma_i-b-a-|f|)\to \bigoplus_{i=1}^sP(-\sigma_i-3) $$are homogeneous maps,  

\item"{(2)}" the map $D_3$ of (3.12) is a lift of $d_3\:\Bbb F_3\to \Bbb F_2$ to $P$, 

\item"{(3)}" both product matrices $D_3D_4$ and $D_4D_3$ are equal to $f$ times an identity matrix, and

\item"{(4)}" the complex
$$\dots \to \Bbb F_3(-|f|)@> \bar D_3>> \Bbb F_2(-|f|) @> \bar D_4 >> \Bbb F_3 @> \bar D_3 >> \Bbb F_2 \to Z\to 0$$ is the minimal homogeneous resolution of $Z$. 
\endroster
\flushpar It follows immediately, that
$$\dots \to \Bbb F_3(-|f|)@> \bar D_3>> \Bbb F_2(-|f|) @> \bar D_4 >> \Bbb F_3 @> \bar D_3 >> \Bbb F_2 @> d_2 >> \Bbb F_1 @> d_1 >> \Bbb F_0 \to R/J\to 0$$ is the minimal homogeneous  resolution of $R/J$. We will modify $D_3$ and $D_4$ to produce the desired matrices $\Phi$ and $\Psi$. 

We next identify an alternating matrix $M$, with entries in $R$,  so that 
$$M\:\bigoplus\limits_{j=1}^{s} R(\sigma_j+3-|f|)@>>> \bigoplus\limits_{i=1}^{s} R(b+a-\sigma_i)\tag 3.13$$ is a homogeneous map
and the maps
$$\bigoplus\limits_{j=1}^{s} R(\sigma_j+3-|f|)@>M>> \bigoplus\limits_{i=1}^{s} R(b+a-\sigma_i)=\Bbb F_2^*@> d_3^* >> \Bbb F_3^*\tag 3.14$$ 
form an exact sequence. 

Fix a basis $y_1,\dots,y_{s}$ for $\Bbb  F_{2}=\bigoplus\limits_{i=1}^{s}R(-(b+a-\sigma_{i}))$, with the degree of $y_i$ equal to $b+a-\sigma_{i}$, and let $y_1^*,\dots,y_{s}^*$ be the corresponding dual basis for $\Bbb  F_{2}^*$. Let $z_1,\dots,z_{s}$ and $\zeta_1,\dots,\zeta_{s}$ be generating sets for $Z$ and $Z^*$, respectively, with $z_i=\operatorname{\frak q}(y_i)$ and $\zeta_i=\ell^{-1}(\pmb s^{-b}(z_i))$ for all $i$. (The shift function is explained in (0.5), $\ell$ is the isomorphism of Proposition 3.10.1, and $\frak q$ is defined in the paragraph before Proposition 3.10.) Form the matrix $M$ with the element  $\zeta_j(z_i)$ of $R$ in position $(i,j)$. Notice that $\zeta_j(z_i)$ has degree $b+2a-\sigma_{i}-\sigma_{j}$; and hence, the map in (3.13) is a homogeneous map. Furthermore, Proposition 3.10.2 shows that  $M$ is an alternating matrix. We have created $M$ so that $M$ carries the $j^{\text{th}}$ basis vector, $v_j$, of $\bigoplus\limits R(\sigma_j+3-|f|)$ to
$\frak q^*(\zeta_j)\in \Bbb F_2^*$. Indeed, $M(v_j)=\sum \zeta_j(z_i) y_i^*\in \Bbb F_2^*$, and, if $y\in \Bbb F_2$, then $M(v_j)$ sends  $y$   to 
$$\align \zeta_j\left(\sum y_i^*(y)\cdot z_i\right)&{}=\zeta_j\left(\sum y_i^*(y)\cdot \frak q(y_i)\right)=
\zeta_j\left(\frak q\left(\sum y_i^*(y)\cdot y_i\right)\right)\\&{}=\zeta_j(\frak q(y))=(\frak q^*(\zeta_j))(y).\endalign$$  
It follows that the image of $M$ in $\Bbb F_2^*$ is equal to the image of $\frak q^*\:Z^*\to \Bbb F^*_2$. On the other hand, we may apply $(\underline{\phantom{X}})^*=\operatorname{Hom}_R(\underline{\phantom{X}}, R)$ to the exact sequence $$\Bbb  F_{3}@> d_{3} >> \Bbb  F_{2}@> \operatorname{\frak q} >> Z\to 0$$ to obtain the exact sequence $$0\to Z^*@> \operatorname{\frak q}^*>> \Bbb  F_{2}^*@> d_{3}^*>> \Bbb  F_{3}^*.$$ We have now shown that (3.14) is an exact sequence. 

We may apply \cite{4, Lemma 2.1} to find a homogeneous invertible matrix $$\varepsilon\: \bigoplus_iP(\sigma_i+3) \to \bigoplus_iP(\sigma_i+3),$$ with  entries in $P$, so that $\varepsilon D_{3}^{\text{\rm T}}$ is an alternating matrix. Define $\Phi$ to be the alternating matrix $D_{3}\varepsilon^{\text{\rm T}}$ and $\Psi$ to be the matrix $(\varepsilon^{\text{\rm T}})^{-1}D_{4}$. It is clear that $(\Phi,\Psi)$ is a matrix factorization of $fI$ and that the complex of (3.2) is a minimal homogeneous resolution. We need only verify that $\Psi$ is an alternating matrix. To do this, we may look   in the quotient field of $P$, where $\Psi$ is equal to $f\Phi^{-1}=(f/\det \Phi)\operatorname{Adj} \Phi$. Observe that $\operatorname{Adj} \Phi$, which is the classical adjoint  of $\Phi$, is an alternating matrix.
 \qed
\enddemo

 \heading Section 4. Examples, further comments, and questions. \endheading
In this section we ask what happens in the situation of Theorem 1.1 when $R$ has small multiplicity and the socle of $R/J$ is pure. Also, we compare the  $R$-resolutions of $R/J$ and $\omega_{R/J}$. 

  In Proposition 4.1, 
 $R$ has small multiplicity and the socle of $R/J$   lives in exactly one degree.
In this case, all of the relevant information (the dimension of the socle of $R/J$, the degrees of the entries of the matrices $d_{3}$ and $d_{4}$, and the degree of the socle elements of $R/J$) is determined by the parity of the back twist $b$ in the $P$-resolution of $P/I$. We give an example for each parity.

\proclaim{Proposition 4.1}
Adopt all of the notation and hypotheses of Theorem  1.1.  Assume  $|f|=3$ and  the socle of $R/J$ lives in exactly one degree.  
 \item{\rm (a)} If $b$ is odd, then $\dim \operatorname{soc} R/J=\frac 32(\mu(I)-1)$, every entry of  $d_{3}$ must have degree $2$,  every entry of $d_{4}$ must have degree $1$,  and $\sigma_{i}=\frac{b-1}2$, for all $i$.
 \item{\rm (b)} If $b$ is even, then $\dim \operatorname{soc} R/J=3(\mu(I)-1)$, every entry of $d_{3}$ must have degree  $1$,  every entry of $d_{4}$ must have degree $2$, and $\sigma_{i}=\frac {b}2-1$, for all $i$. 
\endproclaim

\demo{Proof} Notice that, in the language of Theorem  1.1, $a=a(R)=|f|-3=0$.  Recall that $s=\dim \operatorname{soc} R/J$ and that $Z=\operatorname{im} d_{2}$.  Let $\sigma$ equal to the common value of $\sigma_{i}$,  for $1\le i\le s$, $N_2=b-\sigma$, and $N_3= \sigma+3$.  Theorem 1.1 establishes the exact sequence:
$$0\to Z(-3)\to R(-N_3)^{s}@>d_{3}>> R(-N_2)^{s}@>>> Z\to 0.\tag 4.2$$
Every entry of the matrix $d_{3}$ has the same degree and we denote this degree by $\deg d_{3}$.
It is clear that $$\deg d_{3}=N_3-N_2=2\sigma+3-b.$$In other words, $b+\deg d_{3}-3$ must be even and  $$\sigma=\frac{b+\deg d_{3}-3}2.$$
The matrices $d_{4}$ and $d_{3}$ can be lifted  to $P$ to give  a matrix factorization of   $fI$ and $|f|=3$. So, $\deg d_{3}+\deg d_{4}=3$ and $\deg d_{3}$ is equal to either $1$ or $2$.
If $\deg d_{3}=1$, then  $b$ is even and $\sigma=\frac {b}2-1$. If $\deg d_{3}=2$, then $b$ is odd and $\sigma=\frac{b-1}2$.

Let $e(\underline{\phantom{X}})$ represent multiplicity. When $n$ is large, (4.2) yields
$$\dim Z_n-\dim Z_{n-3}=s\left [\dim R_{n-N_2}-\dim R_{n-N_3}
\right ].\tag 4.3$$The left side of (4.3)
is
$3e(Z)$. 
The right side is 
$$se(R)[N_3-N_2]=se(R)\deg d_{3}.$$The $R$-module $Z$ has positive rank equal to $\mu(I)-1$; therefore, $$e(Z)=e(R)\operatorname{rank}(Z)=e(R)(\mu(I)-1)$$ and  $3(\mu(I)-1)=s\deg d_{3}$.  We have $$s=\frac{3(\mu(I)-1)}{\deg d_{3}},$$with $\deg d_{3}$ equal to $1$ or $2$. Recall that the ideal $I$ is a  grade three Gorenstein ideal; consequently, $\mu(I)$ is automatically odd.
 \qed
\enddemo

\example{Example 4.4} Let $P=k[x,y,z]$, where $k$ is a field of characteristic $p$, $f$ be the polynomial $x^3+y^3+z^3$, $I$ be the grade three Gorenstein ideal $(x^2,xz,xy+z^2,yz,y^2)$ of $P$ (see \cite {2, Proposition 6.2}), $R=P/(f)$ and $J=IR$. The following calculations were made using Macaulay2.
We first give numerical information about the socle and minimal resolution of $R/J^{[p^e]}$, when $p=5$. The notation $2\:\! 5$ under ``pos 1'' next to  $e=0$ means that the module in position $1$ in the minimal $R$-resolution of $R/J$ is $R(-2)^5$. 
The notation $12\:\! 6$ under ``socle'' next to  $e=1$ means that the socle of $R/J^{[p^1]}$ is minimally generated by $6$ generators, each of degree $12$.  The hypotheses of Theorem 1.1 apply to $J^{[p^e]}$ for $1\le e\le 4$. Notice that $b_e=5(5^e)$ is odd, $s_e=\frac 32(\mu(I)-1)$, $\deg d_{e,3}=2$, $\deg_{e,4}=1$, and $\sigma_{e,i}=\frac{b_e-1}2$.
$$\matrix\format \r\ \ &\r\ \ &\r\ \ &\r\ \ &\r\ \ &\r\ \ &\r\ \ &\r\\
e &\text{socle}&\text{pos 0}&\text{pos 1}&\text{pos 2}&\text{pos 3}&\text{pos 4}  \\  
0 &2\:\! 1 &0\:\! 1  &   2\:\! 5 &    3\:\! 6 &   5\:\! 6 &   6\:\! 6       \\
1 &12\:\! 6 &0\:\! 1  &  10\:\! 5 &   13\:\! 6 &  15\:\! 6 &  16\:\! 6 \\
2 &62\:\! 6 &0\:\! 1  &  50\:\! 5 &   63\:\! 6 &  65\:\! 6 &  66\:\! 6 \\
3 &312\:\! 6 &0\:\! 1  & 250\:\! 5 &  313\:\! 6 & 315\:\! 6 & 316\:\! 6 \\
4 &1562\:\! 6 &0\:\! 1  &1250\:\! 5 & 1563\:\! 6 &1565\:\! 6 &1566\:\! 6        \\ 
\endmatrix$$ 
Here is numerical information   when $p=2$. The hypotheses of Theorem 1.1 apply to $J^{[p^e]}$ for $2\le e\le 4$. Notice that $b_e=5(2^e)$ is even, $s_e=3(\mu(I)-1)$, $\deg d_{e,3}=1$, $\deg_{e,4}=2$, and $\sigma_{e,i}=\frac{b_e}2-1 $.
$$\matrix
e&\text{socle}&\text{pos 0}&\text{pos 1}&\text{pos 2}&\text{pos 3}&\text{pos 4}\\
0&2\:\!1&0\:\!1&2\:\!5&3\:\!6&5\:\!6&6\:\!6\\
1&4\:\!7&0\:\!1&4\:\!5&6\:\!12&7\:\!12&9\:\!12\\
2&9\:\!12&0\:\!1&8\:\!5&11\:\!12&12\:\!12&14\:\!12\\
3&19\:\!12&0\:\!1&16\:\!5&21\:\!12&22\:\!12&24\:\!12\\
4&39\:\!12&0\:\!1&32\:\!5&41\:\!12&42\:\!12&44\:\!12\\
\endmatrix$$
\endexample

Return to the situation of Theorem 1.1. We notice that the infinite tails of the $R$-resolutions of $R/J$ and $\omega_{R/J}$ are equal. We wonder how often this phenomenon occurs.
\proclaim{Proposition 4.5} Retain the notation of Theorem 1.1. If $\Bbb L_{\bullet}$ is the minimal $R$-resolution of $\omega_{R/J}$, then the truncation $$\Bbb L_{\ge 3}\:\quad \dots \to \Bbb L_4\to \Bbb L_3$$ is isomorphic, as a complex, to a shift of the truncation $\Bbb  F_{\ge 3}$. \endproclaim
\demo{Proof}    We know the $R$-resolution of $Z$:
$$\dots \to \Bbb  F_{3}\to \Bbb  F_{2}\to Z\to 0.\tag 4.5$$The $R$-resolution of $\operatorname{syz}_2^P(P/I)\otimes_PR$ is given in (3.4).
Consider the short exact sequence (3.6). The map $\alpha_2$ may be lifted to a comparison of complexes from (3.4) to (4.5); the mapping cone of the resulting map of complexes is a resolution of $\omega_{R/J}(-b-a)$:
$$\dots \to \Bbb  F_{5}\to \matrix\Bbb T_3\otimes_PR\\\oplus\\ \Bbb  F_{4}\endmatrix \to \matrix \Bbb T_2\otimes_PR\\\oplus\\ \Bbb  F_{3}\endmatrix  \to 
\Bbb  F_{2}\to \omega_{R/J}(-b-a)\to 0.$$It follows that $\Bbb L_{\ge 3}$ is 
$$\dots \to \Bbb  F_{6}(b+a)\to \Bbb  F_{5}(b+a),$$which is equal to $\Bbb  F_{\ge 3}(-|f|+b+a)=\Bbb  F_{\ge 3}(b-3)$. \qed \enddemo

\medskip

\flushpar {\bf Acknowledgment.}
We appreciate the referee's suggestions; we have followed them and the paper is better because of them.

\enddocument